\documentclass[12pt]{article}
\usepackage{amssymb,latexsym,amsmath}

\begin{document}
\pagestyle{plain} \headheight=5mm  \topmargin=-5mm

\title{A note on Lawson homology for smooth varieties with small Chow groups  }
\author{ Wenchuan Hu }

%\date{}

\maketitle
\newtheorem{Def}{Definition}[section]
\newtheorem{Th}{Theorem}[section]
\newtheorem{Prop}{Proposition}[section]
\newtheorem{Not}{Notation}[section]
\newtheorem{Lemma}{Lemma}[section]
\newtheorem{Rem}{Remark}[section]
\newtheorem{Cor}{Corollary}[section]

\def\s{\section}
\def\ss{\subsection}

\def\d{\begin{Def}}
\def\t{\begin{Th}}
\def\p{\begin{Prop}}
\def\n{\begin{Not}}
\def\la{\begin{Lemma}}
\def\r{\begin{Rem}}
\def\c{\begin{Cor}}
\def\ee{\begin{equation}}
\def\aa{\begin{eqnarray}}
\def\ya{\begin{eqnarray*}}
\def\bd{\begin{description}}

\def\ed{\end{Def}}
\def\et{\end{Th}}
\def\epo{\end{Prop}}
\def\en{\end{Not}}
\def\el{\end{Lemma}}
\def\er{\end{Rem}}
\def\ec{\end{Cor}}
\def\eee{\end{equation}}
\def\eaa{\end{eqnarray}}
\def\ey{\end{eqnarray*}}
\def\ebd{\end{description}}

\def\nn{\nonumber}
\def\bp{{\bf Proof.}\hspace{2mm}}
\def\qe{\hfill$\Box$}
\def\lj{\langle}
\def\rj{\rangle}
\def\dd{\diamond}
\def\ox{\mbox{}}
\def\lb{\label}
\def\rel{\;{\rm rel.}\;}
\def\vp{\varepsilon}
\def\ep{\epsilon}
\def\mod{\;{\rm mod}\;}
\def\exp{{\rm exp}\;}
\def\Lie{{\rm Lie}}
\def\dim{{\rm dim}}
\def\im{{\rm im}\;}
\def\Lag{{\rm Lag}}
\def\Gr{{\rm Gr}}
\def\span{{\rm span}}
\def\Spin{{\rm Spin}}
\def\sign{{\rm sign}\;}
\def\Supp{{\rm Supp}\;}
\def\Sp{{\rm Sp}\;}
\def\ind{{\rm ind}\;}
\def\rank{{\rm rank}\;}
\def\Sg{{\Sp(2n,\C)}}
\def\Na{{\cal N}}
\def\det{{\rm det}\;}
\def\dist{{\rm dist}}
\def\deg{{\rm deg}}
\def\Tr{{\rm Tr}\;}
\def\ker{{\rm ker}\;}
\def\Vect{{\rm Vect}}
\def\H{{\bf H}}
\def\K{{\rm K}}
\def\R{{\bf R}}
\def\C{{\bf C}}
\def\Z{{\bf Z}}
\def\N{{\bf N}}
\def\F{{\bf F}}
\def\Da{{\bf D}}
\def\A{{\bf A}}
\def\La{{\bf L}}
\def\x{{\bf x}}
\def\y{{\bf y}}
\def\Ga{{\cal G}}
\def\Ha{{\cal H}}
\def\L{{\cal L}}
\def\Pa{{\cal P}}
\def\Ua{{\cal U}}
\def\E{{\rm E}}
\def\J{{\cal J}}

\def\m{{\rm m}}
\def\ch{{\rm ch}}
\def\gl{{\rm gl}}
\def\Gl{{\rm Gl}}
\def\Sp{{\rm Sp}}
\def\sf{{\rm sf}}
\def\U{{\rm U}}
\def\O{{\rm O}}
\def\F{{\rm F}}
\def\P{{\rm P}}
\def\D{{\rm D}}
\def\T{{\rm T}}
\def\Sa{{\rm S}}

\begin{center}{\bf \tableofcontents}\end {center}

\begin{center}{\bf Abstract}\end {center}{\hskip .2 in}
Let X be a smooth projective variety of  dimension n on which
rational and homological equivalence coincide for algebraic p-cycles
in the range $0\leq p\leq s$. We show that the homologically trivial
sector of rational Lawson homology $L_pH_k(X,\mathbb{Q})_{hom}$
vanishes for $0\leq n-p\leq s+2$. This is an analogue of a theorem
of C. Peters in  ``dual dimensions". Together with Peters' theorem
we get that the natural transformation
$L_pH_k(X,\mathbb{Q})\rightarrow H_k(X,\mathbb{Q})$ is injective for
all $p$ and $k$ when $X$ is a smooth projective variety of dimension
4 and ${\rm Ch}_0(X)=\Z$.

\s {Introduction} {\hskip .2 in} In this paper, all projective
varieties are defined over $\mathbb{C}$. Let $X$ be a  projective
variety with dimension $n$. Let ${\cal Z}_p(X)$ be the space of
algebraic $p$-cycles on $X$.

The  \textbf{Lawson homology} $L_pH_k(X)$ of $p$-cycles is defined
by

$$
L_pH_k(X) =\left\{
\begin{array}{l}
\pi_{k-2p}({\cal Z}_p(X)), \quad k\geq
2p;\\
0, \quad \quad \quad\quad\quad\quad   k< 2p
\end{array}
\right.
$$
where ${\cal Z}_p(X)$ is given  a natural topology (cf.
\cite{Friedlander}, \cite{Lawson1}). For a general discussion of
Lawson homology, see the survey paper \cite{Lawson2}.

In \cite{Friedlander-Mazur}, Friedlander and Mazur showed that there
are  natural maps, called \textbf{cycle class maps}
 $$ \Phi_{p,k}:L_pH_{k}(X)\rightarrow H_{k}(X). $$

{\Def $L_pH_{k}(X)_{hom}:={\rm
ker}\{\Phi_{p,k}:L_pH_{k}(X)\rightarrow H_{k}(X)\}$;
$L_pH_k(X,\mathbb{Q})_{hom}:=L_pH_{*}(X)_{hom}\otimes {\mathbb{Q}}.$

}

\medskip

C. Peters proved the following result by using the decomposition of
the diagram for the smooth varieties with small Chow groups first
shown by Bloch and Srinivas \cite{Bloch-Srinivas} and generalized by
Paranjape \cite{Paranjape}, Laterveer \cite{Laterveer} and others:

{\Th {\rm (Peters \cite{Peters})} Let $X$ be a smooth projective
variety for which rational and homological equivalence coincide for
$p-$cycles in the range $0\leq p\leq s$ (that is, in the terminology
of [Lat], $X$ has \textbf{small chow groups up to rank} $s$). Then
$L_pH_{*}(X)_{hom}\otimes {\mathbb{Q}}=0$ in the range $0\leq p\leq
s+1$.}

\medskip
By carefully checking the proof of Peters, we discover the symmetry
of the decomposition of the diagonal $\Delta_X\subset X\times X$ and
note that the proof works for $p$-cycles with $0\leq n-p\leq s+2$.

In this note, we will use the tools of Lawson homology  and the
methods and notations given in \cite{Peters} (and the references
therein) to show the following main result:

{\Th Let $X$ be a smooth projective variety of dimension $n$ for
which rational and homological equivalence coincide for $p-$cycles
in the range $0\leq p\leq s$. Then  $L_pH_{*}(X)_{hom}\otimes {
\mathbb{Q}}=0$ in the range $0\leq n-p\leq s+2$.}

\medskip
For convenience, we introduce the following definition:

{\Def A smooth projective variety  $X$ over $\mathbb{C}$ is called
\textbf{rationally connected } if there is  a rational curve through
any 2 points of $X$. A necessary condition for $Z$ to be rationally
connected is that ${\rm Ch}_0(X)\cong \Z$.}

\medskip
For equivalent descriptions of this definition, see the paper of
Koll\'{a}r, Miyaoka and Mori \cite{Kollar-Miyaoka-Mori}.

{\Cor Let $X$ be a smooth projective variety with $\dim (X)=4$ and
${\rm Ch}_0(X)\cong \Z$. Then $L_pH_k(X)_{hom}\otimes
{\mathbb{Q}}=0$ for all $p$ and $k$. In particular, all the smooth
hypersurfaces in $\P^5$ with degree less or equal than 5 have this
property (cf. \cite{Roitman}).}

{\Rem It is shown by the author in [H] that for any smooth
projective rational variety $X$ of \dim$(X)=4$,
$L_pH_k(X)_{hom}=0$ for any $p$ and $k$. Hence the nontriviality
of $L_pH_k(X)_{hom}$ for some $p$,$k$ for a rationally connected
fourfold $X$ would imply irrationality of $X$.}

{\Cor Let $X$ be a general cubic hypersurface of dimension
 less than or equal to 6, then $L_*H_*(X)_{hom}\otimes {\mathbb{Q}}=0$.}

{\Rem Laterveer \cite{Laterveer} showed that Griffiths groups are
torsion for a general cubic hypersurface of dimension less than or
equal to 6. For the general cubic sevenfold in $\P^8$, Albano and
Collino showed that  ${\rm Griff}_3(X)$ (which is $\cong
L_3H_6(X)_{hom}$ by Friedlander in \cite{Friedlander}) is nontrivial
even after tensoring with $\mathbb{Q}$. }

{\Rem This work was done in Spring of 2005 as part of my Ph. D.
thesis. It was included in my research statement and put on my web
page

{\rm \quad  http://www.math.sunysb.edu/\~wenchuan/job/rs.pdf}

\noindent in November 2005. I recently learned that M. Voineagu has
independently obtained this result (cf. \cite{Voineagu}).}

\s {The Proof of the  main result} {\hskip .2 in} The proof of the
main theorem is based on: the Lemma 12 in \cite{Peters}, the
decomposition of the diagonal given in \cite{Paranjape}, and the
computation of Lawson homology of codimension 1 cycles for a smooth
projective variety given by Friedlander \cite{Friedlander}.

For convenience, we write the results we need as follows:

{\Th {\rm (Friedlander \cite{Friedlander})} Let $X$ be any smooth
projective variety of dimension $n$. Then we have the following
isomorphisms
$$
\left\{
\begin{array}{l}
 L_{n-1}H_{2n}(X)\cong \Z,\\
 L_{n-1}H_{2n-1}(X)\cong H_{2n-1}(X,\Z),\\
 L_{n-1}H_{2n-2}(X)\cong H_{n-1,n-1}(X,\Z)=NS(X)\\
L_{n-1}H_{k}(X)=0 \quad for\quad k> 2n.\\

\end{array}
\right.
$$}

{\Rem From this theorem we have $L_{n-1}H_*(X)_{hom}=0$ for any
smooth projective variety $X$ with $\dim(X)=n$.}

\medskip
Now we need to review some definitions about the action of
correspondences. Let $X$ and $Y$ be smooth projective varieties with
$\dim(X)=n$. For $\alpha\in {\cal Z}_{n+d}(X\times Y)$, one puts
 $$\alpha_*(u)=(p_2)_*[p_1^*(u)\cdot \alpha], \quad u\in {\cal Z}_p(X)$$
where $(p_2)_*$ is the proper push-forward, $p_1^*$ is the flat pull
back and  the "." denotes the intersection product of cycles [Pe,
definition 10]. In this way, $\alpha_*$ gives a correspondence
homomorphism
   $$\alpha_*:{\cal Z}_p(X)\rightarrow {\cal Z}_{p+d}(Y).$$
This $\alpha_*$ induces a map (also denoted by $\alpha_*$) on
Lawson homology groups
$$\alpha_*:L_pH_k(X)\rightarrow L_{p+d}H_{k+2d}(Y) $$
which depends only on the class of $\alpha$ in the Chow group of
$X\times Y$ modulo algebraic equivalence. For the details of the
argument here, see [\cite{Peters}, section 1 C.]

 The key Lemma we need was given by Peters as follows:

{\Prop (\cite{Peters}, Lemma 12) Assume that $X$ and $Y$ are smooth
projective varieties and let $\alpha\subset X\times Y$ be  an
irreducible cycle of dimension $\dim(X)=n$, supported on $V\times
W$, where, $V\subset X$ is a subvariety of dimension $v$ and
$W\subset Y$ a subvariety of dimension $w$. Let $\tilde{V}$ , resp.
$\tilde{W}$ be a resolution of singularities of $V$, resp. $W$  and
let $i:\tilde{V}\rightarrow X$ and $j:\tilde{W}\rightarrow Y$ be the
corresponding morphisms. With $\tilde{\alpha}\subset \tilde{V}\times
\tilde{W}$ the proper transform of $\alpha$ and $p_1$, resp. $p_2$
the projections from $X\times Y$ to the first. resp. the second
factor, there is a commutative diagram

$$
\begin{array}{cccc}
  L_{p-n+v+w}H_{k+2(v+w-n)}(\tilde{V}\times \tilde{W}) &
   \stackrel{\tilde{\alpha}_*}{\longrightarrow} & L_pH_k(\tilde{V}\times \tilde{W})   \\

 \uparrow p_1^*&      & \downarrow  (p_2)_*    \\

 L_{p-n+v}H_{k+2(v-n)}(\tilde{V})&   & L_pH_k(\tilde{W})   \\

 \uparrow i^*&   & \downarrow j_*   \\

  L_pH_k(X) & \stackrel{{\alpha}_*}{\longrightarrow} & L_pH_k(Y).

\end{array}
$$
Here $i^*$ is induced by the Gysin homomorphism, $p_1^*$ is the
flat pull-back, and $(p_2)_*$  and $j_*$ come from proper push
forward. In particular, $\alpha_*=0$ if $p<n-v$ or if $p>w$.
Moreover, $\alpha_{n-v}$ acts trivially on $L_{n-v}H_*(X)_{hom}$,
while $\alpha_w$ acts trivially on $L_{w}H_*(X)_{hom}$.
 }

\medskip
There is a corollary of this proposition given by Peters:

{\Cor {\rm (\cite{Peters}, Corollary 13)} An irreducible cycle
$\alpha \subset X\times X$ supported on a product variety $V\times
W$ with $\dim V+ \dim W=n=\dim(X)$ acts trivially on
$L_*H_*(X)_{hom}$.}

\medskip
Combining Friedlander's result (Theorem 2.1) and Peters' Lemma
(Proposition 2.1), we have the following:

{\Cor Under the assumptions of Proposition 2.1, we have that
$\alpha_{w-1}$ acts trivially on $L_{w-1}H_*(X)_{hom}$.}

\medskip
Now we want to recall some results about the decomposition of the
diagonal given in \cite{Bloch-Srinivas} and generalized by Paranjape
\cite{Paranjape} and Laterveer \cite{Laterveer} with more general
triviality hypotheses on the Chow group as stated in Theorem 1.1.
Since the decomposition of diagonal is symmetric, we have the
following version of the diagonal (cf. \cite{Voisin}, Theorem
10.29):

{\Th Let $X$ be a smooth projective variety. Assume that for
$p\leq s$, the maps
$$cl: {\rm CH}_p(X)\otimes {\mathbb{Q}}\rightarrow H^{2n-2p}(X,{\mathbb{Q}})$$
are injective. Then there exists a decomposition
$$\Delta_X= \alpha^{(0)}+\cdots+\alpha^{(s)}+\beta\in {\rm CH}^n(X\times X)\otimes {\mathbb{Q}},$$
where $\alpha^{(p)}$ is supported in $V_p\times W_{n-p}$,
$p=0,\cdots, s$ with $\dim V_p=p$ and $\dim W_{n-p}=n-p$, and
$\beta$ is supported in $X\times W_{n-s-1}$.}

\medskip
Using  the above theorem and Corollary 2.1, we deduce that the
identity acts as $\beta$ on the homologically zero part of the
Lawson homology $L_*H_*(X)_{hom}$. Applying Proposition 2.1 and
Corollary 2.2, we have the following main result:

{\Th Let $X$ be a smooth projective variety such that the maps
$$cl: {\rm CH}_p(X)\otimes {\mathbb{Q}}\rightarrow H^{2n-2p}(X,{\mathbb{Q}})$$
are injective for $p\leq s$. Then $L_{n-p}H_*(X)_{hom}\otimes
{\mathbb{Q}}=0$  for $p=0,\cdots, s+1,s+2$.}

\medskip
As the application, we get Corollary 1.1 immediately.

Recall a result in \cite{Paranjape} and \cite{Schoen}, i.e., the
general cubic hypersurface $X$ of dimension greater than or equal to
5 has ${\rm Ch}_1(X)\cong \Z$ (Certainly ${\rm Ch}_0(X)\cong \Z$ by
Ro\u\i tman \cite{Roitman}.) Hence we have the following

{\Cor Let $X$ be a general cubic hypersurface of dimension
 less than or equal 6, then $L_*H_*(X)_{hom}=0$.}\qe

%\begin{center}{\bf Acknowledge}\end {center} I
%would like to express my gratitude to my advisor, Blaine Lawson,
%for all his help.

\medskip

\noindent
Department of Mathematics,\\
Stony Brook University, SUNY,\\
Stony Brook, NY 11794-3651\\
Email:wenchuan@math.sunysb.edu
\end{document}